\documentclass[10pt,reqno]{amsart}
\usepackage{graphicx}
\usepackage{amsfonts}
\usepackage{amssymb}
\usepackage{amsmath}
\usepackage{amsxtra}
\usepackage{latexsym}
\usepackage{epstopdf}
\usepackage{mathrsfs}
\usepackage{esint}
\usepackage{mathrsfs,galois,booktabs,ctable,threeparttable,tabularx,algorithm,algorithmic}
\usepackage{enumitem}

\newtheorem{theorem}{Theorem}[section]
\newtheorem{lemma}[theorem]{Lemma}

\newtheorem{corollary}[theorem]{Corollary}
\newtheorem{proposition}[theorem]{Proposition}
\newtheorem{Assumption}{Assumption}[section]

\theoremstyle{definition}

\theoremstyle{remark}

\numberwithin{equation}{section}

\begin{document}

\title[]
{Convergence rates of a dual gradient method for
constrained linear ill-posed problems}

\author{Qinian Jin}

\address{Mathematical Sciences Institute, Australian National
University, Canberra, ACT 2601, Australia}
\email{qinian.jin@anu.edu.au} \curraddr{}




\begin{abstract}
In this paper we consider a dual gradient method for solving linear
ill-posed problems $Ax = y$, where $A : X \to  Y$ is a bounded linear operator from a
Banach space $X$ to a Hilbert space $Y$. A strongly convex penalty function is used
in the method to select a solution with desired feature. Under variational source
conditions on the sought solution, convergence rates are derived when the method
is terminated by either an {\it a priori} stopping rule or the discrepancy principle. We
also consider an acceleration of the method as well as its various applications.
\end{abstract}

\def\p{\partial}
\def\l{\langle}
\def\r{\rangle}
\def\C{\mathcal C}
\def\D{\mathscr D}
\def\a{\alpha}
\def\b{\beta}
\def\d{\delta}

\def\la{\lambda}
\def\ep{\varepsilon}
\def\Ga{\Gamma}

\def\bx{\bf x}
\def\by{{\bf y}}
\def\bA{\bf A}
\def\bV{\bf V}
\def\bW{\bf W}
\def\C{\mathcal C}
\def\G{\mathcal G}

\def\N{\mathcal N}
\def\R{\mathcal R}
\def\X{\mathcal X}
\def\Y{\mathcal Y}
\def\B{\mathcal B}
\def\A{\mathcal A}
\def\H{\mathcal H}
\def\bA{{\bf A}}
\def\bx{{\bf x}}
\def\bb{{\bf b}}
\def\ba{{\bf a}}
\def\bV{{\bf V}}
\def\bW{{\bf W}}
\def\bQ{{\bf Q}}
\def\bD{{\bf D}}
\def\D{\mathscr D}
\def\X{{\mathscr X}}
\def\Y{{\mathscr Y}}
\def\R{{\mathcal R}}
\def\E{{\mathcal E}}
\def\EE{{\mathbb E}}
\def\P{{\mathbb P}}

\def\B{\mathcal B}
\def\A{\mathcal A}
\def\a{\alpha}
\def\b{\beta}
\def\d{\delta}
\def\la{\lambda}
\def\R{\mathcal R}
\def\l{\langle}
\def\r{\rangle}
\def\la{\lambda}
\def\p{\partial} 
\def\E{{\mathcal E}}
\maketitle

\section{\bf Introduction}
\setcounter{equation}{0}

Many linear inverse problems can be formulated into the minimization problem
\begin{align}\label{dgm.1}
\min\{\R(x) : x \in X \mbox{ and } Ax = y\}, 
\end{align}
where $A : X \to Y$ is a bounded linear operator from a Banach space $X$ to a Hilbert space $Y$, $y \in \mbox{Ran}(A)$, the range of $A$, and $\R : X \to  (-\infty, \infty]$ is a proper, lower semi-continuous, convex function that is used to select a solution with desired feature. Throughout the paper, all spaces are assumed to be real vector spaces; however, all results still hold for complex vector spaces by minor modifications adapted to complex environments. The norms in $X$ and $Y$ are denoted by the same notation $\|\cdot\|$. We also use the same notation
$\l\cdot, \cdot\r$ to denote the duality pairing in Banach spaces and the inner product in
Hilbert spaces. When the operator $A$ does not have a closed range, the problem (\ref{dgm.1}) is ill-posed in general, thus, if instead of $y$, we only have a noisy data $y^\d$
satisfying
$$
\|y^\d-y\| \le \d
$$
with a small noise level $\d>0$, then replacing $y$ in (\ref{dgm.1}) by $y^\d$ may lead to a problem
that is not well-defined; even if it is well-defined, the solution may not depend
continuously on the data. In order to use a noisy data to find an approximate
solution of (\ref{dgm.1}), a regularization technique should be employed to remove the
instability (\cite{EHN1996,SKHK2012})

In this paper we will consider a dual gradient method to solve (\ref{dgm.1}). This
method is based on applying the gradient method to its dual problem. In order
for a better understanding, we provide a brief derivation of this method which
is well known in optimization community (\cite{B2017,T1991}); the facts from convex analysis
that are used will be reviewed in Section \ref{sect2}. Assume that we only have a noisy
data $y^\d$ and consider the problem (\ref{dgm.1}) with $y$ replaced by $y^\d$. The associated
Lagrangian function is
$$
{\mathcal L}(x, \la) = \R(x) - \l \la, Ax -y^\d\r, \quad x \in  X \mbox{ and } \la \in Y
$$
which induces the dual function
$$
\inf_{x\in X} \left\{\R(x) - \l \la,  Ax -y^\d\r\right\} =  -\R^*(A^*\la) + \l \la, y^\d\r,
$$
where $A^*: Y \to X^*$ denotes the adjoint of $A$ and $\R^*: X^*\to (-\infty, \infty]$ denotes
the Legendre-Fenchel conjugate of $\R$. Thus the corresponding dual problem is
\begin{align}\label{dgm.2}
\min_{\la\in Y} \left\{d_{y^\d}(\la): = \R^*(A^*\la)- \l \la, y^\d\r \right\}. 
\end{align}
Assuming that $\R$ is strongly convex, then $\R^*$ is continuous differentiable with 
$\nabla \R^*: X^*\to X$ and so is the function $\la \to d_{y^\d}(\la)$ on $Y$. 
Therefore, we may apply a gradient method to solve (\ref{dgm.2}) which leads to
$$
\la_{n+1} = \la_n - \gamma (A \nabla \R^*(A^*\la_n)- y^\d) , 
$$
where $\gamma>0$ is a step-size. Let $x_n:= \nabla \R^*(A^*\la_n)$. Then by the properties of
subdifferential we have $A^*\la_n\in \p \R(x_n)$ and hence
$$
x_n \in \arg\min_{x\in X} \left\{\R(x) -\l \la_n, A x- y^\d\r\right\}.
$$
Combining the above two equations results in the following dual gradient method
\begin{align}\label{dgm.3}
\begin{split}
& x_n  = \arg\min_{x\in X} \left\{ \R(x) - \l \la_n, A x- y^\d\r \right\}, \\
& \la_{n+1} = \la_n - \gamma (A x_n - y^\d). 
\end{split}
\end{align}
Note that when $X$ is a Hilbert space and $\R(x) = \|x\|^2/2$, the method (\ref{dgm.3}) 
becomes the standard linear Landweber iteration in Hilbert spaces (\cite{EHN1996}).

By setting $\xi_n := A^*\la_n$, we can obtain from (\ref{dgm.3}) the algorithm
\begin{align}\label{dgm.4}
\begin{split}
& x_n = \arg \min_{x\in X} \left\{\R(x) -\l \xi_n, x\r \right\} , \\
& \xi_{n+1} = \xi_n - \gamma A^*(Ax_n - y^\d) .
\end{split}
\end{align}
Actually the method (\ref{dgm.4}) is equivalent to (\ref{dgm.3}) when the initial guess
$\xi_0$ is chosen from $\mbox{Ran}(A^*)$, the range of $A^*$. Indeed, under the given condition on $\xi_0$, we can conclude from (\ref{dgm.4}) that $\xi_n \in \mbox{Ran}(A^*)$ for all $n$. Assuming $\xi_n = A^* \la_n$ for some $\la_n \in Y$, we can easily see that $x_n$ defined by the first equation in (\ref{dgm.4}) satisfies the first equation in (\ref{dgm.3}). Furthermore, from the second equation in (\ref{dgm.4}) we have 
$$
\xi_{n+1} = A^*(\la_n - \gamma (A x_n - y^\d))
$$
which means $\xi_{n+1} = A^* \la_{n+1}$ with $\la_{n+1}$ defined by the second equation in (\ref{dgm.3}).

The method (\ref{dgm.4}) as well as its generalizations
to linear and nonlinear ill-posed problems in Banach spaces have been considered
in \cite{BH2012,Jin2016,JW2013,KSS2009,RJ2020,SLS2006} and the convergence property has been proved when the
method is terminated by the discrepancy principle. However, except for the linear
and nonlinear Landweber iteration in Hilbert spaces (\cite{EHN1996,HNS1995}), the convergence
rate in general is missing from the existing convergence theory. In this paper we
will consider the dual gradient method (\ref{dgm.3}) and hence the method (\ref{dgm.4}) under
the discrepancy principle
\begin{align}\label{dgm.5}
\| A x_{n_\d} - y^\d\| \le \tau \d < \|A x_n - y^\d\|, \quad 0 \le n< n_\d
\end{align}
with a constant $\tau>1$ and derive the convergence rate when the sought solution
satisfies a variational source condition. This is the main contribution of the present
paper. We also consider accelerating the dual gradient method by Nesterov's acceleration strategy 
and provide a convergence rate result when the method is
terminated by an {\it a priori} stopping rule. Furthermore, we discuss various applications of our convergence theory: 
we provide a rather complete analysis of the
dual projected Landweber iteration for solving linear ill-posed problems in Hilbert
spaces with convex constraint which was proposed in \cite{E1992} with only preliminary results; 
we also propose an entropic dual gradient method using Boltzmann-Shannon
entropy to solve linear ill-posed problems whose solutions are probability density
functions.

In the existing literature there exist a number of regularization methods for solving (\ref{dgm.1}), including the Tikhonov regularization method, the augmented Lagrangian method, and the nonstationary iterated Tikhonov regularization (\cite{FS2010,JS2012,SKHK2012}. In particular, we would like to mention that the augmented Lagrangian method
\begin{align}\label{dgm.6}
\begin{split}
& x_n \in \arg \min_{x\in X} \left\{ \R(x) - \l \la_n, A x- y^\d\r + \frac{\gamma_n}{2} \|Ax - y^\d\|^2 \right\}, \\
& \la_{n+1} = \la_n - \gamma_n (A x_n - y^\d) 
\end{split}
\end{align}
has been considered in \cite{FG2012,FLR2011,FS2010,Jin2017} for solving ill-posed problem (\ref{dgm.1}) as a regularization method. 
This method can be viewed as a modification of the dual gradient
method (\ref{dgm.3}) by adding the augmented term $\frac{\gamma_n}{2} \|A x-y^\d\|^2$ to the definition of $x_n$.
Although the addition of this extra term enables to establish the regularization
property of the augmented Lagrangian method under quite general conditions on
$\R$, it destroys the decomposability structure and thus extra work has to be done
to determine $x_n$ at each iteration step. In contrast, the convergence analysis of the
dual gradient method requires $\R$ to be strongly convex, however the determination of $x_n$ 
is much easier in general. In fact $x_n$ can be given by a closed formula in
many interesting cases; even if $x_n$ does not have a closed formula, there exist fast
algorithms for solving the minimization problem that is used to define $x_n$ since
it does not involve the operator $A$, see Section \ref{sect4} and \cite{BH2012,Jin2016,JW2013} for instance. This
can significantly save the computational time.

The paper is organized as follows, In Section \ref{sect2}, we give a brief review of some
basic facts from convex analysis in Banach spaces. In Section \ref{sect3}, after a quick
account on convergence, we focus on deriving the convergence rates of the dual
gradient method under variational source conditions on the sought solution when
the method is terminated by either an {\it a priori} stopping rule or the discrepancy
principle; we also discuss the acceleration of the method by Nesterov’s strategy.
Finally in Section \ref{sect4}, we address various applications of our convergence theory.

\section{\bf Preliminaries}\label{sect2}
\setcounter{equation}{0}

In this section, we will collect some basic facts on convex analysis in Banach
spaces which will be used in the analysis of the dual gradient method (1.3); for
more details one may refer to \cite{BZ2005,Z2002} for instance.

Let $X$ be a Banach space whose norm is denoted by $\|\cdot\|$, we use $X^*$ to denote its dual space. Given $x\in X$ and $\xi\in X^*$ we write $\l \xi, x\r = \xi(x)$ 
for the duality pairing. For a convex function $f : X \to  (-\infty, \infty]$,  we use
$$
\mbox{dom}(f) := \{x \in X : f(x) < \infty\}
$$
to denote its effective domain. If $\mbox{dom}(f) \ne \emptyset$, $f$ is called proper. 
Given $x\in \mbox{dom}(f)$, an element $\xi\in X^*$ is called a subgradient of $f$ at $x$ if
$$
f(\bar x) \ge  f(x) + \l\xi, \bar x - x\r,  \quad \forall \bar x \in X.
$$
The collection of all subgradients of $f$ at $x$ is denoted as $\p f(x)$ and is called the
subdifferential of $f$ at $x$. If $\p f(x) \ne \emptyset$, then $f$ is called subdifferentiable at $x$. Thus
$x \to \p f(x)$ defines a set-valued mapping $\p f$ whose domain of definition is defined as
$$
\mbox{dom}(\p f) := \{x \in \mbox{dom}(f) : \p f(x) \ne \emptyset\}.
$$
Given $x\in \mbox{dom}(\p f)$ and $\xi \in \p f(x)$, the Bregman distance induced by $f$ at $x$ in
the direction $\xi$ is defined by
$$
D_f^\xi(\bar x, x) := f(\bar x) -  f(x) - \l \xi, \bar x - x\r,  \quad \forall \bar x \in X
$$
which is always nonnegative.

For a proper function $f : X\to (-\infty, \infty]$, its Legendre-Fenchel conjugate is
defined by
$$
f^*(\xi) :=  \sup_{x\in X}  \{\l \xi, x\r - f(x)\}, \quad  \xi \in X^* 
$$
which is a convex function taking values in $(-\infty, \infty]$. According to the definition
we immediately have the Fenchel-Young inequality
\begin{align}\label{dgm.21}
f^*(\xi) + f(x) \ge \l \xi, x\r 
\end{align}
for all $x\in X$ and $\xi\in X^*$. If $f : X\to  (-\infty, \infty]$ is proper, lower semi-continuous
and convex, $f^*$ is also proper and 
\begin{align}\label{dgm.22}
\xi \in \p f(x) \Longleftrightarrow x\in \p f^*(\xi) \Longleftrightarrow  f(x) + f^*(\xi) = \l \xi, x\r.
\end{align}
We will use the following version of the Fenchel-Rockafellar duality formula (see \cite[Theorem 4.4.3]{BZ2005}).

\begin{proposition}\label{dgm.prop21}
Let $X$ and $Y$ be Banach spaces, let $f : X \to (-\infty, \infty]$ and $g : Y \to (-\infty, \infty]$ 
be proper, convex functions, and let $A : X \to Y$ be a bounded linear operator. If there is 
$x_0 \in \emph{dom}(f)$ such that $A x_0 \in \emph{dom}(g)$ and $g$ is continuous at
$A x_0$, then
\begin{align}\label{dgm.FRD}
\inf_{x\in X} \{f(x) + g(Ax)\} = \sup_{\eta \in Y^*} \{-f^*(A^*\eta) - g^*(-\eta) \}.
\end{align}
\end{proposition}

A proper function $f : X \to (-\infty, \infty]$ is called strongly convex if there exists a
constant $\sigma>0$ such that
\begin{align}\label{dgm.23}
f(t\bar x + (1-t) x) + \sigma t(1-t) \|\bar x -x\|^2 \le  tf(\bar x) + (1-t)f(x)
\end{align}
for all $\bar x, x\in \mbox{dom}(f)$ and $t\in [0, 1]$. The largest number $\sigma>0$ such that (\ref{dgm.23})
holds true is called the modulus of convexity of $f$. It can be shown that a proper,
lower semi-continuous, convex function $f : X \to (-\infty, \infty]$ is strongly convex with
modulus of convexity $\sigma>0$ if and only if
\begin{align}\label{dgm.24}
D_f^\xi(\bar x, x) \ge \sigma \|x-\bar x\|^2
\end{align}
for all $\bar x\in \mbox{dom}(f)$, $x\in \mbox{dom}(\p f)$ and $\xi\in \p f(x)$; see \cite[Corollary 3.5.11]{Z2002}.
Furthermore, \cite[Corollary 3.5.11]{Z2002} also contains the following important result
which in particular shows that the strong convexity of $f$ implies the continuous
differentiability of $f^*$.

\begin{proposition}\label{dgm.prop22}
Let $X$ be a Banach space and let $f : X \to (-\infty, \infty]$ be a proper,
lower semi-continuous, strongly convex function with modulus of convexity $\sigma>0$.
Then $\mbox{dom}(f^*) = X^*$, $f^*$ is Fr\'{e}chet differentiable and its gradient 
$\nabla f^*: X^*\to X$ satisfies
$$
\|\nabla f^*(\xi) -\nabla f^*(\eta) \| \le \frac{\|\xi-\eta\|}{2\sigma} 
$$
for all $\xi, \eta \in X^*$. 
\end{proposition} 

It should be emphasized that $X$ in Proposition 2.2 can be an arbitrary Banach
space. For the gradient $\nabla f^*$ of $f^*$, it is in general a mapping from $X^* \to X^{**}$, 
the second dual space of $X$. Proposition 2.2 actually concludes that, for each $\xi\in X^*$, 
$\nabla f^*(\xi)$ is an element in $X^{**}$ that can be identified with an element in $X$ via the
canonical embedding $X \to X^{**}$, and thus $\nabla f^*$ is a mapping from $X^*$ to $X$.

\section{\bf Main results}\label{sect3}
\setcounter{equation}{0}

This section focuses on the study of the dual gradient method (\ref{dgm.3}). We will make
the following assumption.

\begin{Assumption}\label{dgm.ass1}
\begin{enumerate}[leftmargin = 0.8cm]
\item[\emph{(i)}] $X$ is a Banach space, $Y$ is a Hilbert space, and $A : X \to Y$ is a bounded linear operator;

\item[\emph{(ii)}] $\R : X\to (-\infty, \infty]$ is a proper, lower semi-continuous, strongly convex function with
 modulus of convexity $\sigma>0$;

\item[\emph{(iii)}] The equation $Ax = y$ has a solution in $\mbox{dom}(\R)$.
\end{enumerate}
\end{Assumption}

Under Assumption \ref{dgm.ass1}, one can use \cite[Proposition 3.5.8]{Z2002} to conclude that (\ref{dgm.1})
has a unique solution $x^\dag$ and, for each $n$, the minimization problem involved in the
method (\ref{dgm.3}) has a unique minimizer $x_n$ and thus the method (\ref{dgm.3}) is well-defined.
By the definition of $x_n$ we have
\begin{align}\label{dgm.31}
A^* \la_n\in \p \R(x_n). 
\end{align}
By virtue of (\ref{dgm.22}) and Proposition \ref{dgm.prop22} we further have
\begin{align}\label{dgm.32}
x_n = \nabla \R^*(A^*\la_n)
\end{align}
for all $n\ge 0$. 

\subsection{\bf Convergence}\label{sect3.1}

The regularization property of a family of gradient type methods, including (\ref{dgm.4})
as a special case, have been considered in \cite{JW2013} for solving ill-posed problems in
Banach spaces. Adapting the corresponding result to the dual gradient method
(\ref{dgm.3}) we can obtain the following convergence result.

\begin{theorem}\label{dgm.thm31}
Let Assumption \ref{dgm.ass1} hold and let $L := \|A\|^2/(2\sigma)$. Consider the dual
gradient method (\ref{dgm.3}) with $\la_0=0$ for solving (\ref{dgm.1}).

\begin{enumerate}[leftmargin = 0.8cm]
\item[\emph{(i)}]  If $0<\gamma\le 1/L$ then for the integer $n_\d$ chosen such that $n_\d\to \infty$ and
$\d^2 n_\d \to 0$ as $\d\to 0$ there hold
$$
\R(x_{n_\d}) \to \R(x^\dag)  \quad \mbox{ and } \quad D_{\R}^{A^*\la_{n_\d}}(x^\dag, x_{n_\d}) \to 0 
$$
and hence $\|x_{n_\d}-x^\dag\| \to 0$ as $\d \to 0$. 

\item[\emph{(ii)}] If $\tau>1$ and $\gamma>0$ are chosen such that $1-1/\tau-L\gamma>0$, then the discrepancy
principle (\ref{dgm.5}) defines a finite integer $n_\d$ with
$$
\R(x_{n_\d}) \to \R(x^\dag)  \quad \mbox{ and } \quad D_{\R}^{A^*\la_{n_\d}}(x^\dag, x_{n_\d}) \to 0 
$$
and hence$\|x_{n_\d} - x^\dag\|\to 0$ as $\d\to 0$.
\end{enumerate}
\end{theorem}

Theorem \ref{dgm.thm31} gives the convergence results on the method (\ref{dgm.3}) with $\la_0=0$. The convergence result actually holds for any initial guess $\la_0$ with the iterative sequence defined by (\ref{dgm.3}) converging to a solution $x^\dag$ of $A x = y$ with the property 
$$
D_\R^{A^*\la_0}(x^\dag, x_0) = \min\left\{D_\R^{A^*\la_0}(x, x_0): x\in X \mbox{ and } A x = y\right\}, 
$$
where $x_0 =\arg\min_{x\in X} \{\R(x)-\l \la_0,x \r\}$; this can be seen from Theorem \ref{dgm.thm31} 
by replacing $\R(x)$ by $D_\R^{A^*\la_0}(x, x_0)$. This same remark applies to the convergence rate results in the forthcoming subsection. 
For simplicity of exposition, in the following we will consider only the method (\ref{dgm.3}) 
with $\la_0=0$. 

In \cite{JW2013} the convergence result was stated for $X$ to be a reflexive Banach space.
The reflexivity of $X$ was only used in \cite{JW2013} to show the well-definedness of each $x_n$
by the procedure of extracting a weakly convergent subsequence from a bounded
sequence. Under Assumption \ref{dgm.ass1} (ii) the reflexivity of $X$ is unnecessary as the strong
convexity of $\R$ guarantees that each $x_n$ is well-defined in an arbitrary Banach
space, see \cite[Proposition 3.5.8]{Z2002}. This relaxation on $X$ allows the convergence
result to be used in a wider range of applications, see Section \ref{sect4.2} for instance.

The work in \cite{JW2013} actually concentrates on proving part (ii) of Theorem \ref{dgm.thm31}, i.e.
the regularization property of the method terminated by the discrepancy principle
and part (i) was not explicitly stated. However, the argument can be easily adapted
to obtain part (i) of Theorem \ref{dgm.thm31}, i.e. the regularization property of the method
under an a priori stopping rule.

It should be mentioned that the convergence $\R(x_{n_\d}) \to \R(x^\dag)$ was not established in \cite{JW2013}. 
However, if the residual $\|A x_n - y^\d\|$ is monotonically decreasing with respect to $n$, then, 
following the proof in \cite{JW2013}, one can easily establish the convergence $\R(x_{n_\d}) \to \R(x^\dag)$ 
as $\d\to 0$. For the dual gradient method (\ref{dgm.3}), the monotonicity of $\|A x_n - y^\d\|$ is established 
in the following result which is also useful in the forthcoming analysis on deriving convergence rates.

\begin{lemma}\label{dgm.lem32}
Let Assumption \ref{dgm.ass1} hold and let $0<\gamma \le 4 \sigma/\|A\|^2$. Then for the
sequence $\{x_n\}$ defined by (\ref{dgm.3}) there holds
$$
\|A x_{n+1} - y^\d\| \le \|A x_n -y^\d\|
$$
for all integers $n \ge 0$. 
\end{lemma}

\begin{proof}
Recall from (\ref{dgm.31}) that $A^*\la_n\in \p \R(x_n)$ for each $n \ge 0$. By using (\ref{dgm.24}) and
the equation $\la_{n+1} = \la_n - \gamma (A x_n - y^\d)$,  we have
\begin{align*}
2 \sigma \|x_{n+1} - x_n\|^2 
& \le D_{\R}^{A^*\la_n}(x_{n+1}, x_n) + D_{\R}^{A^*\la_{n+1}} (x_n, x_{n+1}) \\
& = \l A^*\la_{n+1} - A^*\la_n, x_{n+1} - x_n\r \\
& = \l \la_{n+1} - \la_n, A x_{n+1} - A x_n\r \\
& = \gamma \l A x_n - y^\d, A x_n - A x_{n+1}\r. 
\end{align*}
In view of the polarization identity in Hilbert spaces, we further have
\begin{align*}
2 \sigma\|x_{n+1} - x_n\|^2 
& \le \frac{\gamma}{2} \left(\|A x_n - y^\d\|^2 - \|A x_{n+1} - y^\d\|^2 + \|A(x_{n+1} - x_n)\|^2\right) \\
& \le  \frac{\gamma}{2} \left(\|A x_n - y^\d\|^2 - \|A x_{n+1} - y^\d\|^2\right)  
+\frac{\gamma \|A\|^2}{2}  \|x_{n+1} - x_n\|^2.
\end{align*}
Since $0<\gamma \le 4 \sigma/\|A\|^2$, we thus obtain the monotonicity of $\|A x_n - y^\d\|^2$ with
respect to $n$. 
\end{proof}

\subsection{\bf Convergence rates}\label{sect3.2}

In this subsection we will derive the convergence rates of the dual gradient method
(\ref{dgm.3}) when the sought solution satisfies certain variational source conditions. The
following result plays a crucial role for achieving this purpose.

\begin{proposition}\label{dgm.prop33}
Let Assumption \ref{dgm.ass1} hold and let $d_{y^\d}(\la):= \R^*(A^*\la) - \l \la, y^\d\r$. 
Let $L := \|A\|^2/(2\sigma)$. Consider the dual gradient method (\ref{dgm.3}) with $\la_0 = 0$. If
$0<\gamma \le 1/L$ then for any $\la\in Y$ there holds
\begin{align*}
d_{y^\d}(\la) - d_{y^\d}(\la_{n+1}) 
& \ge \frac{1}{2\gamma (n+1)} \left(\|\la_{n+1}-\la\|^2 - \|\la\|^2\right) \\
& \quad \, + \left\{\left(\frac{1}{2} - \frac{L\gamma}{4}\right) n + \left(\frac{1}{2} - \frac{L\gamma}{2}\right)\right\} 
\gamma \|A x_n - y^\d\|^2
\end{align*}
for all $n \ge 0$. 
\end{proposition}

\begin{proof}
Since $\R$ is strongly convex with modulus of convexity $\sigma>0$, it follows from Proposition \ref{dgm.prop22} 
that $\R^*$ is continuously differentiable and 
$$
\|\nabla \R^*(\xi) -\nabla \R^*(\eta)\| \le \frac{\|\xi-\eta\|}{2\sigma}, \quad \forall \xi, \eta \in X^*.
$$
Consequently, the function $\la \to d_{y^\d}(\la)$ is differentiable on $Y$ and its gradient is given by 
$$
\nabla d_{y^\d}(\la) = A \nabla \R^*(A^* \la) - y^\d
$$
with 
$$
\|\nabla d_{y^\d}(\tilde \la) - \nabla d_{y^\d}(\la) \| \le L \|\tilde \la -\la\|, \quad \forall \tilde \la, \la \in Y, 
$$
where $L= \|A\|^2/(2\sigma)$. Therefore 
$$
d_{y^\d}(\la_{n+1}) \le d_{y^\d}(\la_n) + \l \nabla d_{y^\d}(\la_n), \la_{n+1} -\la_n\r + \frac{L}{2} \|\la_{n+1} - \la_n\|^2. 
$$
By the convexity of $d_{y^\d}$ we have for any $\la \in Y$ that 
$$
d_{y^\d}(\la_n) \le d_{y^\d}(\la) + \l \nabla d_{y^\d}(\la_n), \la_n -\la\r. 
$$
Combining the above equations we thus obtain
$$
d_{y^\d}(\la_{n+1}) \le d_{y^\d}(\la) + \l \nabla d_{y^\d}(\la_n), \la_{n+1}-\la\r + \frac{L}{2} \|\la_{n+1}-\la_n\|^2.
$$
By using (\ref{dgm.32}) we can see $\nabla d_{y^\d}(\la_n) = A x_n - y^\d$ which together with the equation 
$\la_{n+1} - \la_n = -\gamma (A x_n - y^\d)$ shows that $\nabla d_{y^\d}(\la_n) = (\la_n- \la_{n+1})/\gamma$. 
Consequently 
$$
d_{y^\d}(\la_{n+1}) \le d_{y^\d}(\la) + \frac{1}{\gamma}\l \la_n - \la_{n+1}, \la_{n+1} - \la\r + \frac{L}{2} \|\la_{n+1} - \la_n\|^2. 
$$
Note that 
$$
\l \la_n -\la_{n+1}, \la_{n+1} - \la\r 
= \frac{1}{2} \left(\|\la_n -\la\|^2 -\|\la_{n+1} - \la\|^2 - \|\la_{n+1} - \la_n\|^2\right). 
$$
Therefore 
\begin{align}\label{dgm.33}
d_{y^\d}(\la) - d_{y^\d}(\la_{n+1}) 
& \ge \frac{1}{2\gamma} \left(\|\la_{n+1}-\la\|^2 - \|\la_n-\la\|^2 \right)  \nonumber \\
& \quad \, + \left(\frac{1}{2\gamma}- \frac{L}{2} \right) \|\la_{n+1} - \la_n\|^2. 
\end{align}
Let $m\ge 0$ be any number. By summing (\ref{dgm.33}) over $n$ from $n =0$ to $n = m$ 
and using $\la_0 = 0$ we can obtain
\begin{align}\label{dgm.34}
\sum_{n=0}^m \left(d_{y^\d}(\la)-d_{y^\d}(\la_{n+1})\right) 
& \ge \frac{1}{2\gamma} \left(\|\la_{m+1} -\la\|^2 - \|\la\|^2\right)  \nonumber \\
& \quad \, + \left(\frac{1}{2\gamma} - \frac{L}{2} \right) \sum_{n=0}^m \|\la_{n+1} - \la_n\|^2. 
\end{align}
Next we take $\la = \la_n$ in (\ref{dgm.33}) to obtain 
$$
d_{y^\d}(\la_n) - d_{y^\d}(\la_{n+1}) \ge \left(\frac{1}{\gamma}- \frac{L}{2} \right) \|\la_{n+1} - \la_n\|^2.
$$
Multiplying this inequality by $n$ and then summing over $n$ from $n=0$ to $n = m$ we can obtain
$$
\sum_{n=0}^m n \left(d_{y^\d}(\la_n) - d_{y^\d}(\la_{n+1})\right) 
\ge \left(\frac{1}{\gamma}-\frac{L}{2} \right) \sum_{n=0}^m n \|\la_{n+1} -\la_n\|^2. 
$$
Note that 
$$
\sum_{n=0}^m n \left(d_{y^\d}(\la_n) - d_{y^\d}(\la_{n+1})\right) 
= - (m+1) d_{y^\d}(\la_{m+1}) + \sum_{n=0}^m d_{y^\d}(\la_{n+1}). 
$$
Thus 
$$
- (m+1) d_{y^\d}(\la_{n+1}) + \sum_{n=0}^m d_{y^\d}(\la_{n+1}) 
\ge \left(\frac{1}{\gamma}-\frac{L}{2}\right) \sum_{n=0}^m n \|\la_{n+1} - \la_n\|^2. 
$$
Adding this inequality to (\ref{dgm.34}) gives 
\begin{align*}
(m+1) \left(d_{y^\d}(\la) - d_{y^\d}(\la_{m+1})\right) 
& \ge \frac{1}{2\gamma} \left(\|\la_{m+1}-\la\|^2 - \|\la\|^2\right) \\
& \quad \, + \left(\frac{1}{2\gamma} - \frac{L}{2} \right) \sum_{n=0}^m \|\la_{n+1} - \la_n\|^2 \\
& \quad \, + \left(\frac{1}{\gamma}-\frac{L}{2} \right) \sum_{n=0}^m n \|\la_{n+1} - \la_n\|^2. 
\end{align*}
Recall that $\la_n -\la_{n+1} = \gamma (A x_n- y^\d)$. By using the monotonicity of $\|A x_n - y^\d\|$ 
shown in Lemma \ref{dgm.lem32} we then obtain
\begin{align*}
& (m+1) \left(d_{y^\d}(\la) - d_{y^\d}(\la_{m+1})\right)  \\
& \ge \frac{1}{2\gamma} \left(\|\la_{m+1} - \la\|^2 -\|\la\|^2\right) 
+ \left(\frac{1}{2\gamma} - \frac{L}{2} \right) \gamma^2 \sum_{n=0}^m \|A x_n - y^\d\|^2 \\
& \quad \, + \left(\frac{1}{\gamma}- \frac{L}{2} \right) \gamma^2 \sum_{n=0}^m n \|A x_n - y^\d\|^2 \\
& \ge \frac{1}{2\gamma} \left(\|\la_{m+1}- \la\|^2 -\|\la\|^2 \right) 
+ \left(\frac{1}{2} - \frac{L\gamma}{2}\right) \gamma (m+1) \|A x_m - y^\d\|^2 \\
& \quad \, + \left(1- \frac{L\gamma}{2}\right) \gamma \frac{m(m+1)}{2} \|A x_m - y^\d\|^2. 
\end{align*}
The proof is therefore complete. 
\end{proof}

We now assume that the unique solution $x^\dag$ satisfies a variational source condition specified
in the following assumption. 

\begin{Assumption}\label{dgm.ass2}
For the unique solution $x^\dag$ of (\ref{dgm.1}) there is an error measure function 
$\E^\dag: \mbox{dom}(\R) \to [0, \infty)$ with $\E^\dag(x^\dag) = 0$ such that 
$$
\E^\dag(x) \le \R(x) - \R(x^\dag) + M \|A x - y\|^q, \quad \forall x \in \mbox{dom}(\R)
$$
for some $0<q\le 1$ and some constant $M >0$. 
\end{Assumption}

Variational source conditions were first introduced in \cite{HKPS2007}, as a generalization
of the spectral source conditions in Hilbert spaces, to derive convergence rates of
Tikhonov regularization in Banach spaces. This kind of source conditions was further generalized, 
refined and verified, see \cite{F2018,FG2012,G2010,HM2012,HW2015,HW2017} for instance. The error
measure function $\E^\dag$ in Assumption \ref{dgm.ass2} is used to measure the speed of convergence;
it can be taken in various forms and the usual choice of $\E^\dag$ is the Bregman distance induced by $\R$. 
Use of a general error measure functional has the advantage of covering a wider range of applications. For instance, in reconstructing sparse solutions of ill-posed problems, one may consider the sought solution in the $\ell^1$ space and take $\R(x)=\|x\|_{\ell^1}$.  In this situation,
convergence under the Bregman distance induced by $\R$ may not provide useful approximation result because two points with zero Bregman distance may have arbitrarily large $\ell^1$-distance. However,
under certain natural conditions, the variational source conditions can be verified with $\E^\dag(x)=\|x-x^\dag\|_{\ell^1}$;
see \cite{F2018,FG2018}.

We first derive the convergence rates for the dual gradient method (\ref{dgm.3}) under
an {\it a priori} stopping rule when $x^\dag$ satisfies the variational source condition specified
in Assumption \ref{dgm.ass2}. 

\begin{theorem}\label{dgm.thm34}
Let Assumption \ref{dgm.ass1} hold and let $L: = \|A\|^2/(2\sigma)$. If $0<\gamma \le 1/L$ and $x^\dag$
satisfies the variational source conditions specified in Assumption \ref{dgm.ass2}, then
for the dual gradient method (\ref{dgm.3}) with the initial guess $\la_0= 0$ there holds
$$
\E^\dag(x_n) \le C\left(n^{-\frac{q}{2-q}} + \d^q + n^{\frac{1-q}{2-q}}\d + n \d^2\right)
$$
for all $n \ge 1$, where $C$ is a generic positive constant independent of $n$ and $\delta$.
Consequently, by choosing an integer $n_\d$ with $n_\d\sim \d^{q-2}$ we have
$$
\E^\dag(x_{n_\d}) = O(\d^q).
$$
\end{theorem}

\begin{proof}
Let $d_y(\la):=\R^*(A^*\la)-\l\la, y\r$. Since $0 <\gamma\le 1/L$, from Proposition \ref{dgm.prop33} it
follows that
\begin{align}\label{dgm.35}
 & \left\{\left(\frac{1}{2}-\frac{L\gamma}{4}\right) n +\left(\frac{1}{2}-\frac{L\gamma}{2}\right)\right\} 
\gamma\|A x_n - y^\d\|^2 \nonumber \\
& \le d_{y^\d}(\la) - d_{y^\d}(\la_{n+1}) - \frac{1}{2\gamma (n+1)} \left(\|\la_{n+1}-\la\|^2 -\|\la\|^2\right)  \nonumber \\
& = d_y(\la) - d_y(\la_{n+1}) + \l \la_{n+1}-\la, y^\d-y\r  \nonumber \\
& \quad \, - \frac{1}{2\gamma (n+1)} \left(\|\la_{n+1}-\la\|^2-\|\la\|^2\right)
\end{align}
for all $\la \in Y$. By the Cauchy-Schwarz inequality we have
$$
\l \la_{n+1} -\la, y^\d-y\r \le \d \|\la_{n+1} -\la\| 
\le\frac{1}{4 \gamma (n+1)} \|\la_{n+1}- \la\|^2 + \gamma (n+1) \d^2. 
$$
Thus, it follows from (\ref{dgm.35}) that
\begin{align*}
c_0 n \|A x_n - y^\d\|^2 & + \frac{1}{4\gamma(n + 1)} \|\la_{n+1} - \la\|^2 \\
& \le d_y(\la) - d_y(\la_{n+1})  + \frac{\|\la\|^2}{2\gamma (n + 1)}  +  \gamma (n + 1)\d^2, 
\end{align*}
where $c_0:= (1/2-L\gamma/4) \gamma>0$. By virtue of the inequality $\|\la_{n+1}\|^2 
\le 2 (\|\la\|^2 + \|\la_{n+1}-\la\|^2)$ we then have
\begin{align*}
c_0 n \|A x_n - y^\d\|^2 & + \frac{1}{8\gamma(n + 1)} \|\la_{n+1}\|^2 \\
& \le d_y(\la) - d_y(\la_{n+1})  + \frac{3\|\la\|^2}{4\gamma (n + 1)}  +  \gamma (n + 1)\d^2.
\end{align*}
By the Fenchel-Young inequality (\ref{dgm.21}) and $A x^\dag = y$ we have
\begin{align}\label{dgm.36}
d_y(\la_{n+1}) &= \R^*(A^*\la_{n+1}) - \l\la_{n+1}, A x^\dag\r \nonumber \\
&= \R^*(A^*\la_{n+1}) - \l A^* \la_{n+1}, x^\dag\r \nonumber \\
& \ge - \R(x^\dag) . 
\end{align}
Therefore 
\begin{align*}
& c_0 n \|A x_n - y^\d\|^2  + \frac{1}{8 \gamma (n+1)} \|\la_{n+1}\|^2 \\
& \le \R^*(A^* \la) - \l \la, y\r + \R(x^\dag) + \frac{3\|\la\|^2}{4 \gamma (n+1)} + \gamma (n+1) \d^2
\end{align*}
for all $\la \in Y$. Consequently 
\begin{align*}
& c_0 n \|A x_n - y^\d\|^2  + \frac{1}{8 \gamma (n+1)} \|\la_{n+1}\|^2 \\
& \le \inf_{\la\in Y} \left\{\R^*(A^* \la) - \l \la, y\r + \R(x^\dag) + \frac{3\|\la\|^2}{4 \gamma (n+1)}\right\}  
+ \gamma (n+1) \d^2 \\
& = \R(x^\dag) - \sup_{\la\in Y} \left\{-\R^*(A^* \la) + \l \la, y\r - \frac{3\|\la\|^2}{4 \gamma (n+1)}\right\}  
+ \gamma (n+1) \d^2.
\end{align*}
According to the Fenchel-Rockafellar duality formula given in Proposition \ref{dgm.prop21}, we
have
\begin{align*}
\sup_{\la\in Y} \left\{-\R^*(A^* \la) + \l \la, y\r - \frac{3\|\la\|^2}{4 \gamma (n+1)}\right\}
= \inf_{x\in X} \left\{\R(x) + \frac{1}{3} \gamma (n+1) \|A x- y\|^2\right\}. 
\end{align*}
Indeed, by taking $f(x) = \R(x)$ for $x \in X$ and $g(z) = \frac{1}{3} \gamma (n+1) \|z-y\|^2$ for $z \in Y$, we can obtain this identity immediately from (\ref{dgm.FRD}) by noting that 
$$
g^*(\la) = \frac{3}{4\gamma (n+1)} \|\la\|^2 + \l \la, y\r, 
\quad \la \in Y. 
$$
Therefore 
\begin{align}\label{dgm.37}
& c_0 n \|A x_n - y^\d\|^2  + \frac{1}{8 \gamma (n+1)} \|\la_{n+1}\|^2 
\le \eta_n + \gamma (n+1) \d^2, 
\end{align}
where 
\begin{align}\label{dgm.38}
\eta_n:= \sup_{x\in X} \left\{ \R(x^\dag) - \R(x) -\frac{1}{3} \gamma (n+1) \|A x - y\|^2 \right\}. 
\end{align}
We now estimate $\eta_n$ when $x^\dag$ satisfies the variatioinal source condition given in Assumption \ref{dgm.ass2}. 
By the nonnegativity of $\E^\dag$ we have $\R(x^\dag) - \R(x) \le M \|A x- y\|^q$. Thus
\begin{align}\label{dgm.39}
\eta_n & \le \sup_{x\in X} \left\{ M \|A x - y\|^q - \frac{1}{3} \gamma (n+1) \|A x- y\|^2 \right\} \nonumber \\
& \le \sup_{s\ge 0} \left\{ M s^q - \frac{1}{3} \gamma (n+1) s^2\right\} = c_1 (n+1)^{-\frac{q}{2-q}},
\end{align}
where $c_1:= \left(1-\frac{q}{2}\right) \left(\frac{3qM}{2\gamma}\right)^{\frac{q}{2-q}} M>0$. 
Combining this with (\ref{dgm.37}) gives 
\begin{align*}
& c_0 n \|A x_n - y^\d\|^2  + \frac{1}{8 \gamma (n+1)} \|\la_{n+1}\|^2 
\le c_1 (n+1)^{-\frac{q}{2-q}} + \gamma (n+1) \d^2
\end{align*}
which implies that 
\begin{align}\label{dgm.310}
\|A x_n - y^\d\| \le C \left(n^{-\frac{1}{2-q}} + \d\right) \quad \mbox{ and } \quad 
\|\la_n\| \le C \left(n^{\frac{1-q}{2-q}} + n \d\right).
\end{align}
Recall $A^* \la_n \in \p \R(x_n)$ from (\ref{dgm.31}), we have 
$$
\R(x_n)- \R(x^\dag) \le \l A^* \la_n, x_n -x^\dag\r = \l \la_n, A x_n - y\r. 
$$
Therefore, by using the variational source condition specified in Assumption \ref{dgm.ass2}, we obtain
\begin{align*}
\E^\dag(x_n)  & \le \R(x_n) -\R(x^\dag) + M \|A x_n -y\|^q \\
& \le \l \la_n, A x_n - y\r + M \|A x_n - y\|^q \\
& \le \|\la_n\|\|A x_n - y\| + M \|A x_n - y\|^q. 
\end{align*}
Thus, it follows from (\ref{dgm.310}) that
\begin{align*}
\E^\dag(x_n) &\le \|\la_n\|\left(\|A x_n - y^\d\| + \d\right) + M \left(\|A x_n - y^\d\| + \d\right)^q \\
& \le C \left(n^{\frac{1-q}{2-q}} + n\d\right) \left(n^{-\frac{1}{2-q}} + \d\right) + C \left(n^{-\frac{1}{2-q}} + \d\right)^q\\
& \le C \left(n^{-\frac{q}{2-q}} + \d^q + n^{\frac{1-q}{2-q}} \d + n \d^2\right). 
\end{align*}
The proof is thus complete. 
\end{proof}

During the proof of Theorem \ref{dgm.thm34}, we have introduced the quantity $\eta_n$ defined by (\ref{dgm.38}).
Taking $x = x^\dag$ in (\ref{dgm.38}) shows $\eta_n \ge 0$. As can be seen from the proof of Theorem \ref{dgm.thm34}, we have 
$$
\eta_n = \inf_{\la\in Y} \left\{ \R^*(A^* \la) - \l \la, y \r + \R(x^\dag) + \frac{3\|\la\|^2}{4\gamma (n+1)}\right\}
$$
by the Fenchel-Rockafellar duality formula. Taking $\la=0$ in this equation gives $0\le \eta_n \le \R(x^\dag) + \R^*(0)<\infty$. The proof of Theorem \ref{dgm.34} demonstrates that $\eta_n$ can decay to $0$ at certain rate if $x^\dag$ satisfies a variational source condition.  

\begin{corollary}\label{dgm.cor35}
Let Assumption \ref{dgm.ass1} hold and let $L := \|A\|^2/(2\sigma)$. If $0<\gamma \le 1/L$
and if there is $\la^\dag \in Y$ such that $A^*\la^\dag \in \p \R(x^\dag)$, then for the dual gradient method
(\ref{dgm.3}) with the initial guess $\la_0 =0$ there holds
\begin{align}\label{dgm.311}
D_{\R}^{A^*\la^\dag} (x_n, x^\dag) \le C \left(n^{-1} + \d + n \d^2\right) 
\end{align}
for all $n \ge 1$, where $C$ is a generic positive constant independent of $n$ and $\d$.
Consequently, by choosing an integer $n_\d$ with $n_\d\sim \d^{-1}$ we have
\begin{align}\label{dgm.312}
D_{\R}^{A^*\la^\dag} (x_{n_\d}, x^\dag) = O(\d) 
\end{align}
and hence $\|x_{n_\d} - x^\dag\| = O(\d^{1/2})$.
\end{corollary}

\begin{proof}
We show that $x^\dag$ satisfies the variational source condition specified in Assumption \ref{dgm.ass2} 
with $q = 1$. The argument is well-known, see \cite{HKPS2007} for instance. Since
$A^*\la^\dag \in \p \R(x^\dag)$ for some $\la^\dag \in Y$, we have for all $x \in \mbox{dom}(\R)$ that
\begin{align*}
D_{\R}^{A^*\la^\dag} (x, x^\dag) 
& = \R(x) - \R(x^\dag) -\l \la^\dag, A x - y\r \\
& \le \R(x) - \R(x^\dag) + \|\la^\dag\| \|A x- y\|
\end{align*}
which shows that Assumption \ref{dgm.ass2} holds with $\E^\dag(x) = D_{\R}^{A^*\la^\dag} (x, x^\dag)$, 
$M = \|\la^\dag\|$ and $q = 1$. Thus by invoking Theorem \ref{dgm.thm34}, we immediately obtain (\ref{dgm.311}) 
which together with the choice $n_\d\sim \d^{-1}$ implies (\ref{dgm.312}). By using (\ref{dgm.24}) we then obtain
$\|x_{n_\d} - x^\dag\| = O(\d^{1/2})$. 
\end{proof}

We next turn to deriving convergence rates of the dual gradient method (\ref{dgm.3})
under the variational source condition given in Assumption \ref{dgm.ass2} when the method
is terminated by the discrepancy principle (\ref{dgm.5}). We will use the following consequence 
of Proposition \ref{dgm.prop33}.

\begin{lemma}\label{dgm.lem36}
Let Assumption \ref{dgm.ass1} hold and let $L := \|A\|^2/(2\sigma)$. Consider the dual
gradient method (\ref{dgm.3}) with $\la_0=0$. If $\tau>1$ and $\gamma>0$ are chosen such that
$1-1/\tau^2- L\gamma>0$, then there is a constant $c_2>0$ such that
$$
c_2 (n + 1)\d^2\le \eta_n \quad \mbox{ and } \quad \frac{1}{8\gamma (n+1)} \|\la_{n+1}\|^2 
\le \eta_n + \gamma (n + 1) \d^2
$$
for all integers $0\le n \le n_\d$, where $n_\d$ is the integer determined by the discrepancy
principle (\ref{dgm.5}) and $\eta_n$ is the quantity defined by (\ref{dgm.38}).
\end{lemma}

\begin{proof}
The second estimate follows directly from (\ref{dgm.37}), actually it holds for all integers $n\ge 0$. It remains only to show the first estimate. For any $n<n_\d$ we have $\|A x_n - y^\d\|>\tau \d$. Therefore from (\ref{dgm.35}) it follows for all $\la \in Y$ that
\begin{align*}
& \left\{\left(\frac{1}{2}-\frac{L\gamma}{4}\right) n + \left(\frac{1}{2}-\frac{L\gamma}{2}\right)\right\} \gamma \tau^2 \d^2 \\
& \le d_y(\la) - d_y(\la_{n+1}) +\l\la_{n+1}-\la, y^\d - y\r - \frac{1}{2\gamma (n+1)} \left(\|\la_{n+1}-\la\|^2 - \|\la\|^2\right).
\end{align*}
By the Cauchy-Schwarz inequality we have
\begin{align*}
& \l \la_{n+1}-\la, y^\d- y\r \\
& \le \d \|\la_{n+1}-\la\| \le \frac{1}{2\gamma (n+1)} \|\la_{n+1}-\la\|^2 + \frac{1}{2} \gamma (n+1) \d^2.
\end{align*}
Therefore 
\begin{align*}
& \left[\left\{\left(\frac{1}{2}-\frac{L\gamma}{4}\right) n + \left(\frac{1}{2}-\frac{L\gamma}{2}\right)\right\} \tau^2 -\frac{1}{2}(n+1)\right]  \gamma \d^2 \\
& \le d_y(\la) - d_y(\la_{n+1}) + \frac{1}{2\gamma (n+1)} \|\la\|^2. 
\end{align*}
By the conditions on $\gamma$ and $\tau$, it is easy to see that
$$
\left[\left\{\left(\frac{1}{2}-\frac{L\gamma}{4}\right) n + \left(\frac{1}{2}-\frac{L\gamma}{2}\right)\right\} \tau^2 -\frac{1}{2}(n+1)\right]\gamma \ge c_2 (n+1), 
$$
where $c_2:= \left((1/2-L\gamma/2)\tau^2-1/2\right) \gamma >0$. Therefore 
$$
c_2(n+1) \d^2 \le d_y(\la) - d_y(\la_{n+1}) + \frac{1}{2\gamma (n+1)} \|\la\|^2. 
$$
According to (\ref{dgm.36}) we have $d_y(\la_{n+1}) \ge -\R(x^\dag)$. Thus
$$
c_2(n+1) \d^2 \le \R^*(A^*\la) -\l \la, y\r + \R(x^\dag) + \frac{1}{2\gamma(n+1)} \|\la\|^2
$$
which is valid for all $\la\in Y$. Consequently
\begin{align*}
c_2(n+1) \d^2 & \le \inf_{\la\in Y} \left\{ \R^*(A^*\la) -\l \la, y\r + \R(x^\dag) + \frac{1}{2\gamma(n+1)} \|\la\|^2\right\} \\
& = \R(x^\dag) - \sup_{\la\in Y} \left\{-\R^*(A^*\la) +\l \la, y\r - \frac{1}{2\gamma(n+1)} \|\la\|^2\right\}. 
\end{align*}
According to the Fenchel-Rockafellar duality formula given in Proposition \ref{dgm.prop21}, we can further obtain
\begin{align*}
c_2 (n+1) \d^2 
& \le \R(x^\dag) - \inf_{x\in X} \left\{\R(x) + \frac{1}{2} \gamma (n+1) \|A x- y\|^2\right\} \\
& = \sup_{x\in X} \left\{\R(x^\dag) - \R(x) - \frac{1}{2} \gamma (n+1) \|A x- y\|^2\right\}\\
& \le \eta_n
\end{align*}
which shows the first estimate. 
\end{proof}

Now we are ready to show the convergence rate result for the dual gradient
method (\ref{dgm.3}) under Assumption \ref{dgm.ass2} when the method is terminated 
by the discrepancy principle (\ref{dgm.5}).

\begin{theorem}\label{dgm.thm37}
Let Assumption \ref{dgm.ass1} hold and let $L := \|A\|^2/(2\sigma)$. Consider the dual
gradient method (\ref{dgm.3}) with the initial guess $\la_0= 0$. Assume that $\tau>1$ and $\gamma>0$ 
are chosen such that $1- 1/\tau^2 - L \gamma>0$ and let $n_\d$ be the integer determined by the
discrepancy principle (\ref{dgm.5}). If $x^\dag$ satisfies the variational source condition specified
in Assumption \ref{dgm.ass2}, then
\begin{align}\label{dgm.313}
\E^\dag(x_{n_\d}) = O(\d^q). 
\end{align}
Consequently, if there is $\la^\dag \in Y$ such that $A^* \la^\dag \in \p \R(x^\dag)$, then
\begin{align}\label{dgm.314}
D_{\R}^{A^*\la^\dag} (x_{n_\d}, x^\dag) = O(\d) 
\end{align}
and hence $\|x_{n_\d} - x^\dag \| = O(\d^{1/2})$.
\end{theorem}

\begin{proof}
By using the variational source condition on $x^\dag$ 
specified in Assumption \ref{dgm.ass2},
the convexity of $\R$, and the fact $A^*\la_{n_\d} \in \p \R(x_{n_\d})$ we have
\begin{align*}
\E^\dag(x_{n_\d}) & \le \R(x_{n_\d}) - \R(x^\dag) + M\|A x_{n_\d} - y\|^q \\
& \le\l\la_{n_\d}, A x_{n_\d} - y\r + M \|A x_{n_\d}-y\|^q \\
& \le \|\la_{n_\d}\| \|A x_{n_\d} - y\| + M \| A x_{n_\d}-y\|^q.  
\end{align*}
By the definition of $n_\d$ we have $\|A x_{n_\d} - y^\d\| \le \tau \d$ and thus
$$
\|A x_{n_\d}-y\| \le \|A x_{n_\d}-y^\d\| + \|y^\d - y\| \le (\tau+1) \d.
$$
Therefore 
\begin{align}\label{dgm.315}
\E^\dag(x_{n_\d}) \le (\tau+1) \|\la_{n_\d}\| \d + M (\tau +1)^q \d^q. 
\end{align}
If $n_\d = 0$, then we have $\la_{n_\d} = 0$ and hence $\E^\dag(x_{n_\d}) \le M (\tau+1)^q \d^q$. In the following
we consider the case $n_\d\ge 1$. We will use Lemma \ref{dgm.lem36} to estimate $\|\la_{n_\d}\|$. By virtue
of Assumption \ref{dgm.ass2} we have $\eta_n \le c_1 (n + 1)^{-\frac{q}{2-q}}$, see (\ref{dgm.39}). Combining this with the estimates in Lemma \ref{dgm.lem36} we can obtain
\begin{align}\label{dgm.316}
c_2 (n+1)^{\frac{2}{2-q}} \d^2 \le c_1    
\end{align}
and 
\begin{align}\label{dgm.317}
\|\la_{n+1}\|^2 \le 8 \gamma c_1 (n+1)^{\frac{2(1-q)}{2-q}} + 8 \gamma (n+1)^2 \d^2 
\end{align}
for all $0 \le n <n_\d$. Taking $n = n_\d- 1$ in (\ref{dgm.316}) gives
$$
n_\d \le \left(\frac{c_1}{c_2 \d^2}\right)^{\frac{2-q}{2}}
$$
which together with (\ref{dgm.317}) with $n = n_\d-1$ shows that
$$
\|\la_{n_\d}\| \le c_3 \d^{q-1}, 
$$
where $c_3 := \sqrt{8\gamma (1+ c_2)(c_1/c_2)^{2-q}}$. Combining this estimate with (\ref{dgm.315}) we finally
obtain
$$
\E^\dag(x_{n_\d}) \le \left(c_3(\tau+1) + M(\tau+1)^q\right) \d^q
$$
which shows (\ref{dgm.313}).

When $A^* \la^\dag \in \p \R(x^\dag)$ for some $\la^\dag \in Y$, we know from the proof of Corollary \ref{dgm.cor35} that Assumption \ref{dgm.ass2} is satisfied with $\E^\dag(x) = D_\R^{A^*\la^\dag}(x, x^\dag)$ and $q = 1$. Thus, we
may use (\ref{dgm.313}) to conclude (\ref{dgm.314}). 
\end{proof}

\subsection{Acceleration}\label{sect3.3}

The dual gradient method, which generalizes the linear Landweber iteration in
Hilbert spaces, is a slowly convergent method in general. To make it more practical important, it is necessary to consider accelerating this method with faster
convergence speed. Since the dual gradient method is obtained by applying the
gradient method to the dual problem, one may consider to accelerate this method
by applying any available acceleration strategy for gradient methods among which
Nesterov's acceleration strategy (\cite{N1983,BT2009,AP2016}) is the most prominent. By applying Nesterov's 
accelerated gradient method to minimize the function $d_{y^\d}(\la) = \R^*(A^*\la) - \l \la, y^\d\r$ 
it leads to the iteration scheme
\begin{align*}
& \hat \la_n = \la_n + \frac{n-1}{n+\a} (\la_n-\la_{n-1}), \\
& \la_{n+1} = \hat \la_n - \gamma \nabla d_{y^\d}(\hat \la_n).  
\end{align*}
Let $\hat x_n = \nabla \R^*(A^* \hat \la_n)$. Then $\nabla d_{y^\d}(\hat \la_n) = A \hat x_n - y^\d$ and $A^*\hat \la_n \in \p \R(\hat x_n)$ which imply that 
\begin{align} \label{dgm.318}
\begin{split}
& \hat \la_n = \la_n + \frac{n-1}{n+\a} (\la_n -\la_{n-1}),\\
& \hat x_n = \arg\min_{x\in X} \left\{ \R(x) - \l \hat \la_n, A x- y^\d\r \right\}, \\
& \la_{n+1} = \hat \la_n - \gamma (A \hat x_n - y^\d), \\
& x_{n+1} = \arg\min_{x\in X} \left\{\R(x) - \l \la_{n+1}, A x - y^\d\r \right\},
\end{split}
\end{align}
where $\la_{-1} = \la_0=0$, $\a\ge 2$ is a given number, and $\gamma>0$ is a step size. We
have the following convergence rate result when the method is terminated by an {\it a priori} stopping rule.

\begin{theorem}\label{dgm.thm38}
Let Assumption \ref{dgm.ass1} hold and let $L := \|A\|^2/(2\sigma)$. Consider the accelerated dual gradient method (\ref{dgm.318}) with noisy data $y^\d$ satisfying $\|y^\d-y\| \le \d$. Assume that $0 <\gamma \le 1/L$ and $\a\ge 2$. If $x^\dag$ satisfies the source condition
$A^*\la^\dag \in \p \R(x^\dag)$ for some $\la^\dag \in Y$, then there exist positive constants $c_4$ and $c_5$
depending only on $\gamma$ and $\a$ such that
\begin{align} \label{dgm.319}
D_{\R}^{A^*\la_n} (x^\dag, x_n)\le \left( \frac{c_4 \|\la^\dag\|}{n} + c_5 n\d\right)^2
\end{align}
for all $n \ge 1$. Consequently by choosing an integer $n_\d$ with $n_\d\sim \d^{-1/2}$ we have
$$
D_{\R}^{A^*\la_{n_\d}} (x^\dag, x_{n_\d}) = O(\d)
$$
and hence $\|x_{n_\d} - x^\dag \| = O(\d^{1/2})$ as $\d\to 0$.
\end{theorem}

\begin{proof}
According to the definition of $x_n$ we have $A^*\la_n\in \p \R(x_n)$ for all $n \ge 1$.
From this fact and the condition $A^*\la^\dag \in \p \R(x^\dag)$ it follows from (\ref{dgm.22}) that
\begin{align}\label{dgm.320}
D_\R^{A^*\la_n}(x^\dag, x_n) 
& = \R(x^\dag) - \R(x_n) - \l \la_n, y -A x_n\r \nonumber \\
& = \left\{\l A^* \la^\dag, x^\dag\r - \R^*(A^* \la^\dag)\right\} 
- \left\{\l A^* \la_n, x_n\r - \R^*(A^*\la_n)\right\} \nonumber\\
& \quad \, - \l \la_n, y - A x_n\r \nonumber\\
& = \R^*(A^*\la_n) - \R^*(A^*\la^\dag) - \l \la_n, y\r + \l \la^\dag, y\r \nonumber \\
& = d_y(\la_n) - d_y(\la^\dag),
\end{align}
where $d_y(\la) := \R^*(A^*\la) - \l \la, y\r$. We need to estimate $d_y(\la_n) - d_y(\la^\dag)$. This can
be done by using a perturbation analysis of the accelerated gradient method, see \cite{AP2016,AD2015}. For completeness, we include a derivation here. Because $A^*\la^\dag \in \p \R(x^\dag)$, we have $x^\dag = \nabla \R^*(A^* \la^\dag)$. Thus
$$
\nabla d_y(\la^\dag) = A \nabla \R^*(A^* \la^\dag) - y = A x^\dag - y = 0.
$$
Since $d_y$ is convex, this shows that $\la^\dag$ is a global minimizer of $d_y$ over $Y$. Note that $\nabla d_{y^\d}(\la) = \nabla d_y(\la) + y- y^\d$. Thus, it follows from the definition of $\la_{n+1}$ that
$$
\la_{n+1} = \hat \la_n - \gamma \left(\nabla d_y(\hat \la_n) + y - y^\d\right). 
$$
Based on this, the Lipschitz continuity of $\nabla d_y$ and the convexity of $\R$, we may
use a similar argument in the proof of Proposition \ref{dgm.prop33} to obtain for any $\la\in Y$ that
\begin{align*}
d_y(\la_{n+1}) & \le d_y(\la) + \l \nabla d_y(\hat \la_n), \la_{n+1}-\la\r + \frac{L}{2} \|\la_{n+1} - \hat \la_n\|^2 \\
& = d_y(\la) + \left\l \frac{1}{\gamma} (\hat \la_n - \la_{n+1}) - (y - y^\d), \la_{n+1} - \la\right\r + \frac{L}{2} \|\la_{n+1} - \hat \la_n\|^2 \\
& = d_y(\la) + \frac{1}{2\gamma} \left(\|\hat \la_n-\la\|^2 -\|\la_{n+1}-\la\|^2 \right)-\l y-y^\d, \la_{n+1} - \la\r \\
& \quad \, - \left(\frac{1}{2\gamma} - \frac{L}{2} \right) \|\la_{n+1} -\hat \la_n\|^2.
\end{align*}
Since $0<\gamma \le 1/L$ and $\|y^\d - y\| \le \d$, we have 
\begin{align}\label{dgm.321}
d_y(\la_{n+1}) \le d_y(\la) + \frac{1}{2\gamma} \left(\|\hat \la_n-\la\|^2 -\|\la_{n+1}-\la\|^2 \right)+ \d\|\la_{n+1} - \la\|.
\end{align}
Note that $\frac{n-1}{n+\a} = \frac{t_n-1}{t_{n+1}}$ with $t_n = \frac{n+\a-1}{\a}$. Now we take $\la = \left(1-\frac{1}{t_{n+1}}\right) \la_n + \frac{1}{t_{n+1}} \la^\dag$ in (\ref{dgm.321}) and use the convexity of $d_y$ to obtain
\begin{align*}
d_y(\la_{n+1}) & \le \left(1-\frac{1}{t_{n+1}}\right) d_y(\la_n) + \frac{1}{t_{n+1}} d_y(\la^\dag) \\
& \quad \, + \frac{1}{2\gamma t_{n+1}^2} \left\|\la^\dag - \left(\la_n + t_{n+1}(\hat \la_n - \la_n)\right)\right\|^2 \\
& \quad \, - \frac{1}{2 \gamma t_{n+1}^2} \left\|\la^\dag - \left(\la_n + t_{n+1}(\la_{n+1}-\la_n)\right)\right\|^2 \\
& \quad \, + \frac{\d}{t_{n+1}} \left\| \la^\dag - \left(\la_n + t_{n+1} (\la_{n+1} - \la_n)\right) \right\|. 
\end{align*}
Let $u_n = \la_{n-1} + t_n (\la_n - \la_{n-1})$. Then it follows from $\hat \la_n = \la_n + \frac{t_n-1}{t_{n+1}} (\la_n -\la_{n-1})$ that $\la_n + t_{n+1} (\hat \la_n - \la_n) = u_n$. Therefore
\begin{align*}
d_y(\la_{n+1}) & \le \left(1-\frac{1}{t_{n+1}}\right) d_y(\la_n) + \frac{1}{t_{n+1}} d_y(\la^\dag) + \frac{1}{2\gamma t_{n+1}^2} \|\la^\dag - u_n\|^2 \\
& \quad \, - \frac{1}{2\gamma t_{n+1}^2} \|\la^\dag - u_{n+1}\|^2 + \frac{\d}{t_{n+1}} \|\la^\dag - u_{n+1}\|. 
\end{align*}
Multiplying both sides by $2\gamma t_{n+1}^2$, regrouping the terms and setting $w_n := d_y(\la_n) - d_y(\la^\dag)$, we obtain
\begin{align}\label{dgm.322}
2 \gamma t_{n+1}^2 w_{n+1} - 2 \gamma t_n^2 w_n 
& \le 2 \gamma \rho_n w_n + \|\la^\dag - u_n\|^2 - \|\la^\dag - u_{n+1}\|^2 \nonumber \\
& \quad \, + 2 \gamma \d t_{n+1} \|\la^\dag - u_{n+1}\|
\end{align}
for all $n \ge 0$, where $\rho_n:= t_{n+1}^2 - t_{n+1} - t_n^2$. Note that $\a\ge 2$ implies $\rho_n \le 0$ for
$n \ge 1$. Let $m \ge 1$ be any integer. Summing the above inequality over $n$ from $n = 1$ to $n = m - 1$ and using $w_n \ge 0$, we can obtain
\begin{align*}
2\gamma t_m^2 w_m + \|\la^\dag - u_m\|^2 
\le 2\gamma t_1^2 w_1 + \|\la^\dag - u_1\|^2 
+ 2 \gamma \d \sum_{k=2}^m t_k \|\la^\dag -u_k\|. 
\end{align*}
Using (\ref{dgm.322}) with $n = 0$ and noting that $t_0^2 +\rho_0 =0$ and $u_0= 0$ we can further obtain
\begin{align}\label{dgm.323}
& 2\gamma t_m^2 w_m + \|\la^\dag - u_m\|^2 \nonumber \\
& \le 2 \gamma (t_0^2 + \rho_0) w_0 + \|\la^\dag\|^2 
+ 2 \gamma \d t_1 \|\la^\dag - u_1\| 
+ 2 \gamma \d \sum_{k=2}^m t_k \|\la^\dag - u_k\| \nonumber\\
& = \|\la^\dag\|^2 + 2 \gamma \d \sum_{k=1}^m t_k \|\la^\dag - u_k\|. 
\end{align}
According to (\ref{dgm.323}), we have
\begin{align}\label{dgm.324}
\|\la^\dag - u_m\|^2 
\le  \|\la^\dag\|^2 + 2 \gamma \d \sum_{k=1}^m t_k \|\la^\dag - u_k\|
\end{align}
from which we may use an induction argument to obtain
\begin{align}\label{dgm.325}
\|\la^\dag - u_m\| 
\le  \|\la^\dag\| + 2 \gamma \d \sum_{k=1}^m t_k
\end{align}
for all integers $m \ge 0$. Indeed, since $u_0 = 0$, (\ref{dgm.325}) holds trivially for $m = 0$.
Assume next that (\ref{dgm.325}) holds for all $0 \le m \le n$ for some $n \ge 0$. We show (\ref{dgm.325}) also holds for $m = n+1$. If there is $0 \le m \le n$ such that $\|\la^\dag - u_{n+1}\| \le  \|\la^\dag - u_m\|$, then by the induction hypothesis we have
$$
\|\la^\dag - u_{n+1}\| \le \|\la^\dag\| + 2\gamma \d \sum_{k=1}^m t_k \le \|\la^\dag\| + 2 \gamma \d \sum_{k=1}^{n+1} t_k. 
$$
So we may assume $\|\la^\dag - u_{n+1}\| > \|\la^\dag - u_m\|$ for all $0 \le m \le n$. It then follows
from (\ref{dgm.324}) that
$$
\|\la^\dag - u_{n+1}\|^2 \le \|\la^\dag\|^2 + 2\gamma \d \left(\sum_{k=1}^{n+1} t_k\right) \|\la^\dag - u_{n+1}\|. 
$$
By using the elementary inequality ``$a^2 \le b^2 + ca \Longrightarrow a \le b+c$ for $a, b, c \ge 0$", we obtain again
$$
\|\la^\dag - u_{n+1}\| \le \|\la^\dag\| + 2\gamma \d \sum_{k=1}^{n+1} t_k. 
$$
By the induction principle, we thus obtain (\ref{dgm.325}). Based on (\ref{dgm.323}) and (\ref{dgm.325}) we have
\begin{align*}
2 \gamma t_m^2 w_m & \le \|\la^\dag\|^2 + 2 \gamma \d \sum_{k=1}^m t_k \|\la^\dag - u_k\| \\
& \le \|\la^\dag\|^2 + 2 \gamma \d \left(\sum_{k=1}^m t_k\right)\left(\|\la^\dag\| + 2 \gamma \d \sum_{k=1}^m t_k\right) \\
& \le \left(\|\la^\dag\|+ 2 \gamma \d \sum_{k=1}^m t_k\right)^2. 
\end{align*}
Thus, by the definition of $t_n$ it is straightforward to see that
$$
d_y(\la_m) - d_y(\la^\dag) 
\le \frac{1}{2\gamma t_m^2} \left(\|\la^\dag\| + 2\gamma \d \sum_{k=1}^m t_k\right)^2 
\le \left(\frac{c_4 \|\la^\dag\|}{m} + c_5 m \d\right)^2, 
$$
where $c_4$ and $c_5$ are two positive constants depending only on $\gamma$ and $\a$. Combining
this with (\ref{dgm.320}) we thus complete the proof of (\ref{dgm.319}). 
\end{proof}

From Theorem \ref{dgm.thm38} it follows that, under the source condition $A^*\la^\dag \in \p \R(x^\dag)$,
we can obtain the convergence rate $\|x_{n_\d} - x^\dag\|
= O(\d^{1/2})$ within $O(\d^{-1/2})$ iterations for the method (\ref{dgm.318}). For the dual gradient method (\ref{dgm.3}), however, we need to perform $O(\d^{-1})$ iterations to achieve the same convergence rate, see Corollary \ref{dgm.cor35}. This demonstrates that the method (\ref{dgm.318}) indeed has acceleration effect.

We remark that Nesterov’s acceleration strategy was first proposed in \cite{Jin2016}
to accelerate gradient type regularization method for linear as well as nonlinear ill-posed problems in Banach spaces and various numerical results were reported which demonstrate the striking performance; see also \cite{HR2017,HR2018,K2021,N2017,ZWJ2019} for
further numerical simulations. Although we have proved in Theorem \ref{dgm.thm38} a convergence rate result for the method (\ref{dgm.318}) under an {\it a priori} stopping rule, the
regularization property of the method under the discrepancy principle is not yet
established for general strongly convex $\R$. However, when $X$ is a Hilbert space
and $\R(x) = \|x\|^2/2$, the regularization property of the corresponding method has
been established in \cite{K2021,N2017} based on a general acceleration framework in \cite{H1991} using
orthogonal polynomials; in particular it was observed in \cite{K2021} that the parameter $\a$
plays an interesting role in deriving order optimal convergence rates. For an analysis of Nesterov’s acceleration for nonlinear ill-posed problems in Hilbert spaces,
one may refer to \cite{HR2018}.

\section{\bf Applications}\label{sect4}
\setcounter{equation}{0}

Various applications of the dual gradient method (\ref{dgm.3}), or equivalently the method
(\ref{dgm.4}), have been considered in \cite{BH2012,Jin2016,JW2013} for sparsity recovery 
and image reconstruction through the choices of $\R$ as strong convex perturbations of 
the $L^1$ and the total variation functionals and the numerical results demonstrate its nice 
performance. In the following we will provide some additional applications.

\subsection{Dual projected Landweber iteration}\label{sect4.1}

We first consider the application of our convergence theory to linear ill-posed
problems in Hilbert spaces with convex constraint. Such problems arise from a
number of real applications including the computed tomography \cite{N2001} in which the
sought solutions are nonnegative.

Let $A : X \to Y$ be a bounded linear operator between two Hilbert spaces $X$
and $Y$ and let $C\subset X$ be a closed convex set. Given $y\in Y$ and assuming that
$Ax = y$ has a solution in $C$, we consider finding the unique solution $x^\dag$ of $Ax = y$
in $C$ with minimal norm which can be stated as the minimization problem
\begin{align}\label{dgm.41}
\min\left\{\frac{1}{2} \|x\|^2 : x \in C \mbox{ and } Ax = y\right\}.
\end{align}
This problem takes the form (\ref{dgm.1}) with $\R(x) := \frac{1}{2} \|x\|^2 + \d_C(x)$, where $\d_C$ denotes
the indicator function of $C$, i.e. $\d_C(x) = 0$  if $x\in C$ and $\infty$ otherwise. Clearly $\R$
satisfies Assumption \ref{dgm.ass1} (ii). It is easy to see that for any $\xi\in X$ the unique solution
of
$$
\min_{x\in C} \left\{\R(x) -\l \xi, x\r\right\} 
$$
is given by $P_C(\xi)$, where $P_C$ denotes the metric projection of $X$ onto $C$. 
Therefore, applying the algorithm (\ref{dgm.3}) to (\ref{dgm.41}) leads to the dual projected Landweber iteration
\begin{align}\label{dgm.42}
\begin{split}
x_n &= P_C (A^*\la_n), \\
\la_{n+1} &  = \la_n - \gamma (A x_n - y^\d)
\end{split}
\end{align}
that has been considered in \cite{E1992}. Besides a stability estimate, it has been shown
in \cite{E1992} that, for the method (\ref{dgm.42}) with exact data $y^\d=y$, if $x^\dag \in P_C (A^*Y)$ then
$$
\sum_{n=1}^\infty \|x_n - x^\dag\|^2 <\infty
$$
which implies $\|x_n - x^\dag\|\to 0$ as $n \to \infty$ but does not provide an error estimate
unless $\|x_n -x^\dag\|$ is monotonically decreasing which is unfortunately unknown.
Therefore, the work in \cite{E1992} does not tell much information about the regularization
property of the method (\ref{dgm.42}). It is natural to ask if the method (\ref{dgm.42}) renders a
regularization method under {\it a priori} or {\it a posteriori} stopping rules and if it is
possible to derive error estimates under suitable source conditions on the sought
solution. Applying our convergence theory can provide satisfactory answers to
these questions with a rather complete analysis on the method (\ref{dgm.42}), see Corollary
\ref{dgm.cor41} below, which is far beyond the one provided in \cite{E1992}. In particular we can obtain
the convergence and convergence rates when the method (\ref{dgm.42}) is terminated by
either an {\it a priori} stopping rule or the discrepancy principle (\ref{dgm.5}).

\begin{corollary}\label{dgm.cor41}
For the linear ill-posed problem (\ref{dgm.41}) in Hilbert spaces constrained
by a closed convex set $C$, consider the dual projected Landweber iteration (\ref{dgm.42}) with
$\la_0 =0$ and with noisy data $y^\d$ 
satisfying $\|y^\d - y\|\le \d$.

\begin{enumerate}[leftmargin = 0.8cm]
\item[\emph{(i)}]  If $0<\gamma \le 1/\|A\|^2$ then for the integer $n_\d$ satisfying $n_\d\to \infty$ and $\d^2 n_\d\to 0$ as $\d\to 0$ there holds $\|x_{n_\d}-x^\dag\| \to 0$ as $\d\to 0$. If in addition $x^\dag$ satisfies the
projected source condition
\begin{align}\label{dgm.43}
x^\dag = P_C((A^*A)^{\nu/2} \omega) \mbox{ for some } 0 < \nu \le 1 \mbox{ and } \omega \in X, 
\end{align}
then with the choice $n_\d\sim \d^{-\frac{2}{1+\nu}}$ we have $\|x_{n_\d} - x^\dag \| = O(\d^{\frac{\nu}{1+\nu}})$.

\item[\emph{(ii)}] If $\tau>1$ and $\gamma>0$ are chosen such that $1 - 1/\tau - \|A\|^2 \gamma >0$, then the
discrepancy principle (\ref{dgm.5}) defines a finite integer $n_\d$ with $\|x_{n_\d} - x^\dag \| \to 0$ as $\d\to 0$. If in addition $x^\dag$ satisfies the projected source condition (\ref{dgm.43}), then
$\|x_{n_\d} - x^\dag \| = O(\d^{\frac{\nu}{1+\nu}})$.
\end{enumerate}
\end{corollary}

\begin{proof}
According to Theorem \ref{dgm.thm31}, Theorem \ref{dgm.thm34} and Theorem \ref{dgm.thm37}, it remains only
to show that, under the projected source condition (\ref{dgm.43}), $x^\dag$ satisfies the variational
source condition
\begin{align}\label{dgm.44}
\frac{1}{4} \|x-x^\dag\|^2 \le \R(x) - \R(x^\dag) + c_\nu \|\omega\|^{\frac{2}{1+\nu}} \|A x - y\|^{\frac{2\nu}{1+\nu}}
\end{align}
for all $x \in \mbox{dom}(\R) = C$, where $c_\nu := 2^{-\frac{2\nu}{1+\nu}} (1+\nu) (1-\nu)^{\frac{1-\nu}{1+\nu}}$. To see this, note
first that for any $x \in C$ there holds
\begin{align*}
\frac{1}{2} \|x-x^\dag\|^2 - \R(x) + \R(x^\dag) 
& = \frac{1}{2} \|x-x^\dag\|^2 - \frac{1}{2} \|x\|^2 + \frac{1}{2} \|x^\dag\|^2 \\
& = \l x^\dag, x^\dag -x\r. 
\end{align*}
By using $x^\dag = P_C ((A^* A)^{\nu/2} \omega)$ and the property of the projection $P_C$ we have
$$
\l (A^*A)^{\nu/2}\omega - x^\dag, x - x^\dag \r \le 0, \quad \forall x \in C 
$$
which implies that 
$$
\l x^\dag, x^\dag -x\r \le \l (A^*A)^{\nu/2} \omega, x^\dag - x\r = \l \omega, (A^*A)^{\nu/2} (x^\dag-x)\r.
$$
Therefore 
\begin{align*}
\frac{1}{2} \|x-x^\dag\|^2 - \R(x) + \R(x^\dag) 
& \le \l \omega, (A^*A)^{\nu/2} (x^\dag-x)\r \\
& \le \|\omega\| \|(A^*A)^{\nu/2} (x^\dag-x)\|.
\end{align*}
By invoking the interpolation inequality (\cite{EHN1996}) and the Young’s inequality we can further obtain
\begin{align*}
\frac{1}{2} \|x-x^\dag\|^2 - \R(x) + \R(x^\dag) 
& \le \|\omega\| \|x-x^\dag\|^{1-\nu} \|A(x-x^\dag)\|^\nu \\
& \le c_\nu \|\omega\|^{\frac{2}{1+\nu}} \|A x - y\|^{\frac{2\nu}{1+\nu}} + \frac{1}{4} \|x-x^\dag\|^2
\end{align*}
which shows (\ref{dgm.44}). The proof is therefore complete. 
\end{proof}

We remark that the projected source condition (\ref{dgm.43}) with $\nu=1$, i.e. $x^\dag \in P_C(A^* Y)$, was first used in \cite{N1988} to derive the convergence rate of Tikhonov regularization in Hilbert spaces with convex constraint.

The dual projected Landweber iteration (\ref{dgm.42}) can be accelerated by Nesterov’s acceleration strategy. As was derived in Section \ref{sect3.3}, the accelerated scheme takes the form
\begin{align}\label{dgm.45}
\begin{split}
& \hat \la_n = \la_n + \frac{n-1}{n+\a} (\la_n -\la_{n-1}), \qquad \hat x_n = P_C(A^* \hat \la_n), \\
& \la_{n+1} = \hat \la_n - \gamma (A \hat x_n - y^\d), 
\qquad x_{n+1} = P_C (A^* \la_{n+1}).
\end{split}
\end{align}
By noting that $\p \R(x) = x + \p \d_C(x)$, it is easy to see that an element $\la^\dag\in Y$ is such that $A^* \la^\dag \in \p \R(x^\dag)$ if and only if $x^\dag = P_C(A^*\la^\dag)$. 
Therefore, by using Theorem \ref{dgm.thm38}, we can obtain the following convergence rate result of the method
(\ref{dgm.45}).

\begin{corollary} \label{dgm.cor42}
For the problem (\ref{dgm.41}) in Hilbert spaces constrained by a closed convex set $C$, consider the method (\ref{dgm.45}) with $\la_0 = \la_{-1} = 0$. If $0 <\gamma \le 1/\|A\|^2$, $\a\ge 2$ and $x^\dag \in P_C(A^* Y)$, then with the choice $n_\d \sim \d^{-1/2}$ we have
$$
\|x_{n_\d} - x^\dag \| =  O(\d^{1/2}) 
$$
as $\d \to 0$. 
\end{corollary}

\subsection{An entropic dual gradient method}\label{sect4.2}

Let $\Omega \subset {\mathbb R}^d$ be a bounded domain and let $A : L^1(\Omega) \to Y$ be a bounded linear
operator, where $Y$ is a Hilbert space. For an element $y \in Y$ in the range of $A$, we consider the equation $Ax = y$. We assume that the sought solution $x^\dag$ is a
probability density function, i.e. $x^\dag \ge 0$ a.e. on $\Omega$ and $\int_\Omega x^\dag = 1$. We may find
such a solution by considering the convex minimization problem
\begin{align}\label{dgm.46}
\min\left\{\R(x) : = f(x) + \d_\Delta (x): x\in L^1(\Omega) \mbox{ and } A x = y\right\},  
\end{align}
where $\d_\Delta$ denotes the indicator function of the closed convex set
\begin{align*}
\Delta := \left\{x\in L^1(\Omega): x\ge 0 \mbox{ a.e. on } \Omega \mbox{ and } \int_\Omega x =1\right\} 
\end{align*}
in $L^1(\Omega)$ and $f$ denotes the negative of the Boltzmann-Shannon entropy, i.e.
$$
f(x) := \left\{\begin{array}{lll}
\int_\Omega x \log x & \mbox{ if } x \in L_+^1(\Omega) \mbox{ and } x \log x \in L^1(\Omega), \\[1.2ex]
\infty & \mbox{ otherwise }
\end{array}\right.
$$
where, here and below, $L_+^p(\Omega):= \{x \in L^p(\Omega): x \ge 0 \mbox{  a.e. on } \Omega \}$ for each $1\le p\le \infty$. The Boltzmann-Shannon entropy has been used in Tikhonov regularization as a stable functional to determine nonnegative solutions; see \cite{AH1991,E1993,EL1993,JH1996} for instance.

In the following we summarize some useful properties of the negative of the Boltzmann-Shannon entropy $f$:

\begin{enumerate}[leftmargin = 0.8cm]
\item[(i)] $f$ is proper, lower semi-continuous and convex on $L^1(\Omega)$; see \cite{AH1991,E1993}.

\item[(ii)] $f$ is subdifferentiable at $x\in L^1(\Omega)$  if and only if $x \in L_+^\infty(\Omega)$ and is bounded
away from zero, i.e.
$$
\mbox{dom}(\p f) = \{x \in L_+^\infty(\Omega): x\ge \beta \mbox{ on } \Omega \mbox{ for some constant } \beta>0\}.
$$
Moreover for each $x\in \mbox{dom}(\p f)$ there holds $\p f(x) = \{1 + \log x\}$; see \cite[Proposition 2.53]{BP2012}.

\item[(iii)] By straightforward calculation one can see that for any $x \in \mbox{dom}(\p f)$ and $\tilde x\in \mbox{dom}(f)$, the Bregman distance induced by $f$ is the Kullback-Leibler
functional
$$
D(\tilde x, x) := \int_\Omega 
\left(\tilde x \log \frac{\tilde x}{x} - \tilde x + x\right).
$$

\item[(iv)] For any $x \in \mbox{dom}(\p f)$ and $\tilde x \in \mbox{dom}(f)$ there holds (see \cite[Lemma 2.2]{BL1991})
\begin{align}\label{dgm.412}
\|x-\tilde x\|_{L^1(\Omega)}^2 \le \left(\frac{4}{3} \|x\|_{L^1(\Omega)} + \frac{2}{3} \|\tilde x\|_{L^1(\Omega)}\right) D(\tilde x, x).
\end{align}
\end{enumerate}

Based on these facts, we can see that the function $\R$ defined in (\ref{dgm.46}) satisfies
Assumption \ref{dgm.ass1}. In order to apply the dual gradient method (\ref{dgm.3}) to solve (\ref{dgm.46}), we need to determine the closed formula for the solution of the minimization problem
involved in the algorithm. By the Karush-Kuhn-Tucker theory, it is easy to see
that, for any $\ell\in L^\infty(\Omega)$, the unique minimizer of
$$
\min\left\{\int_\Omega(x\log x - \ell x): x \ge 0 \mbox{ a.e. on } \Omega \mbox{ and } \int_\Omega x = 1\right\}
$$
is given by $\hat x = e^\ell/\int_\Omega e^\ell$. Therefore we can obtain from the algorithm (\ref{dgm.3}) the following entropic dual gradient method
\begin{align}\label{dgm.47}
\begin{split}
x_n &= \frac{1}{\int_\Omega e^{A^*\la_n}} e^{A^*\la_n}, \\
\la_{n+1} & = \la_n - \gamma (A x_n - y^\d).
\end{split}
\end{align}
We have the following convergence result.

\begin{corollary}\label{dgm.cor43} 
For the convex problem (\ref{dgm.46}), consider the entropic dual gradient method (\ref{dgm.47}) with $\la_0 = 0$ and with noisy data $y^\d$ satisfying $\|y^\d - y\|\le \d$.

\begin{enumerate}[leftmargin = 0.8cm]
\item[\emph{(i)}] If $0 <\gamma \le 1/\|A\|^2$ then for the integer $n_\d$ satisfying $n_\d\to \infty$ and $\d^2 n_\d\to 0$ as $\d\to 0$ there holds $\|x_{n_\d} - x^\dag \|\to 0$ as $\d\to 0$. If in addition $x^\dag$ satisfies the source condition
\begin{align}\label{dgm.48}
1 + \log x^\dag = A^*\la^\dag \quad \mbox{ for some } \la^\dag \in Y, 
\end{align}
then with the choice $n_\d\sim \d^{-1}$ we have $\|x_{n_\d} -x^\dag\|_{L^1(\Omega)} = O(\d^{1/2})$.

\item[\emph{(ii)}] If $\tau>1$ and $\gamma>0$ are chosen such that $1-2/\tau-\|A\|^2 \gamma >0$, then the
discrepancy principle (\ref{dgm.5}) defines a finite integer $n_\d$ with $\|x_{n_\d} - x^\dag \|_{L^1(\Omega)} \to 0$ as $\d\to 0$. If in addition $x^\dag$ satisfies the source condition (\ref{dgm.48}), then $\|x_{n_\d} - x^\dag\|_{L^1(\Omega)}  = O(\d^{1/2})$.
\end{enumerate}
\end{corollary}

\begin{proof}
Under (\ref{dgm.48}) there holds $A^* \la^\dag \in \p f(x^\dag)$. Therefore, by using (\ref{dgm.412}), we have for any $x\in \mbox{dom}(\R)$ that
\begin{align*}
\frac{1}{2} \|x-x^\dag\|_{L^1(\Omega)}^2 
& \le D(x, x^\dag) = f(x) - f(x^\dag) - \l A^* \la^\dag, x-x^\dag\r \\
& = \R(x)-\R(x^\dag) - \l \la^\dag, A x - y\r \\
& \le \R(x) - \R(x^\dag) + \|\la^\dag\| \|A x - y\|,
\end{align*}
where we used $\R(x) = f(x) $ and $\int_\Omega x =1$ for $x \in \mbox{dom}(\R)$. Thus $x^\dag$ satisfies the varaitional source condition specified in Assumption \ref{dgm.ass2} with $\E^\dag(x) = \frac{1}{2} \|x-x^\dag\|_{L^1(\Omega)}^2$, $M =\|\la^\dag\|$ and $q =1$. Now we can complete the proof by
applying Theorem \ref{dgm.thm31}, Corollary \ref{dgm.cor35} and Theorem \ref{dgm.thm37} to the method (\ref{dgm.47}).
\end{proof}

The source condition (\ref{dgm.48}) has been used in \cite{E1993,EL1993} in which one may find
further discussions. We would like to mention that an entropic Landweber method of the form
\begin{align}\label{dgm.49}
x_{n+1} = \frac{x_n e^{\gamma A^*(y^\d-A x_n)}}{\int_\Omega x_n e^{\gamma A^*(y^\d - A x_n)}}
\end{align}
has been proposed and studied in the recent paper \cite{BRB2020} in which weak convergence
in $L^1(\Omega)$ is proved without relying on source conditions and, under the source
condition (\ref{dgm.48}), an error estimate is derived when the method is terminated by an
{\it a priori} stopping rule. Our method (\ref{dgm.47}) is different from (\ref{dgm.49}) due to its primal-dual nature. As stated in Corollary \ref{dgm.cor43}, our method (\ref{dgm.47}) enjoys nicer convergence
properties: it admits strong convergence in $L^1(\Omega)$ in general and, when the source
condition (\ref{dgm.48}) is satisfied, an error estimate can be derived when the method is
terminated by either an {\it a priori} stopping rule or the discrepancy principle.

Applying Nesterov’s acceleration strategy, we can accelerate the entropic dual gradient method (\ref{dgm.47}) by the following scheme
\begin{align}\label{dgm.410}
\begin{split}
& \hat \la_n = \la_n + \frac{n-1}{n+\a} (\la_n - \la_{n-1}), \qquad \hat x_n = \frac{1}{\int_\Omega e^{A^*\hat \la_n}} e^{A^*\hat \la_n}, \\
& \la_{n+1} = \hat \la_n - \gamma (A \hat x_n - y^\d), 
\qquad x_{n+1} = \frac{1}{\int_\Omega e^{A^*\la_{n+1}}} e^{A^*\la_{n+1}}.
\end{split}
\end{align}
By using Theorem \ref{dgm.thm38} we can obtain the following convergence rate result on the method (\ref{dgm.410}) with noisy data.

\begin{corollary} \label{dgm.cor44}
For the minimization problem (\ref{dgm.46}), consider the method (\ref{dgm.410}) with $\la_0 = \la_{-1} = 0$. 
If $0 <\gamma \le 1/\|A\|^2$, $\a\ge 2$ and $x^\dag$ satisfies the source condition (\ref{dgm.48}), then with the choice $n_\d\sim \d^{-1/2}$ we have
$$
\|x_{n_\d} - x^\dag \|_{L^1(\Omega)} = O(\d^{1/2})
$$
as $\d\to 0$. 
\end{corollary}

\section{\bf Conclusion}

Due to its simplicity and relatively small complexity per iteration, Landweber
iteration has received extensive attention in the inverse problems community. In
recent years, Landweber iteration has been extended to solve inverse problems
in Banach spaces with general uniformly convex regularization terms and various
convergence properties have been established. However, except for the linear and
nonlinear Landweber iteration in Hilbert spaces, the convergence rate in general
is missing from the existing convergence theory.

This paper attempts to fill in this gap by providing a novel technique to derive
convergence rates for a class of Landweber type methods. We considered a class of
ill-posed problems defined by a bounded linear operator from a Banach space to
a Hilbert space and used a strongly convex regularization functional to select the
sought solution. The dual problem turns out to have a smooth objective function
and thus can be solved by the usual gradient method. The resulting method is
called a dual gradient method which is a special case of the Landweber type
method in Banach spaces. Applying gradient methods to the dual problem allows
us to interpret the method in a new perspective which enables us to use tools
from convex analysis and optimization to carry out the analysis. We have actually
obtained the convergence and convergence rates of the dual gradient method when
it is terminated by either an {\it a priori} stopping rule or the discrepancy principle.
Furthermore, by applying Nesetrov’s acceleration strategy to the dual problem we
proposed an accelerated dual gradient method and established a convergence rate
result under an {\it a priori} stopping rule. We also discussed some applications, in
particular, as a direct application of our convergence theory, we provided a rather
complete analysis of the dual projected Landweber iteration of Eicke for which
only preliminary result is available in the existing literature.

There are a few of questions which might be interesting for future development.

\begin{enumerate}[leftmargin = 0.8cm]
\item[(i)]  We established convergence rate results for the dual gradient method (\ref{dgm.3})
which require $A$ to be a bounded linear operator and $Y$ a Hilbert space.
Is it possible to establish a general convergence rate result for Landweber
iteration for solving linear as well as nonlinear ill-posed problems in Banach
spaces?

\item[(ii)] For the dual gradient method (\ref{dgm.3}), its analysis under the {\it a priori} stopping
rule allows to take the step-size as $0<\gamma \le 1/L$, while the analysis under the
discrepancy principle (\ref{dgm.5}) requires $\tau>1$ and $\gamma>0$ to satisfy $1-1/\tau^2-L\gamma >0$
which means either $\tau$ has to be large or $\gamma$ has to be small. Is it possible to
develop a convergence theory of the discrepancy principle under merely the
conditions $\tau>1$ and $0<\gamma \le 1/L$?

\item[(iii)] In Section \ref{sect3.3} we considered the accelerated dual gradient method (3.18) and
established a convergence rate result under an {\it a priori} stopping rule. Is it
possible to establish the convergence and convergence rate result of (3.18)
under the discrepancy principle?
\end{enumerate}

\section*{\bf Acknowledgements}

The work of Q Jin is partially supported by the Future Fellowship of the Australian
Research Council (FT170100231).

\bibliographystyle{amsplain}

\end{document}